\documentclass[12pt]{article}
\usepackage{latexsym}
\usepackage{amsfonts}
\usepackage[all]{xy}

\def\Bbb{\mathbb}

\newtheorem{defn}{Definition}
\newtheorem{thm}{Theorem}[section]

\newtheorem{cor}[thm]{Corollary}
\newtheorem{lem}[thm]{Lemma}

\newtheorem{question}[thm]{Question}
\newenvironment{proof}{\medskip \noindent
{\bf Proof.}}{\hfill \rule{.5em}{1em}
\\}

\newenvironment{rmk}{\mbox{ }\\{\bf  Remark}\mbox{ }}{
\hfill $\Box$\mbox{}\bigskip}
\begin{document}
\sloppy
\title{Surgery, Yamabe invariant, \\ and Seiberg-Witten theory}

\author{Chanyoung Sung\thanks{email address: cysung@kias.re.kr \
 }
\\
National Institute for Mathematical Sciences\\
385-16 Doryong-dong Yuseong-gu Daejeon  Korea}
\date{}
\maketitle

\begin{abstract}
By using the gluing formula of the Seiberg-Witten invariant, we
compute the Yamabe invariant $Y(X)$ of $4$-manifolds $X$ obtained by
performing surgeries along points, circles or tori on compact K\"ahler
surfaces. For instance, if $M$ is a compact K\"ahler
surface of nonnegative Kodaira dimension, and $N$ is a smooth closed oriented $4$-manifold with $b_2^+(N)=0$ and $Y(N)\geq 0$, then we show that
$$Y(M\# N)=Y(M).$$

\end{abstract}

\section{Introduction}
The Yamabe invariant is a real-valued invariant of a smooth closed manifold
defined using the scalar curvature. It somehow measures how much
the negative scalar curvature is inevitable, and it can be used as
a means to get to a canonical metric on a given manifold.

 Let $M$ be a closed smooth $n$-manifold. In any conformal class
$$[g]=\{\varphi g \mid \varphi:M \rightarrow \Bbb
R^+ \  \textrm{is smooth} \},$$ there exists  a smooth Riemannian
metric of constant scalar curvature, so-called \emph{Yamabe
metric}, realizing the minimum of the normalized total scalar
curvature
$$\inf_{\tilde{g} \in
[g]} \frac{\int_M s_{\tilde{g}}\ dV_{\tilde{g}}}{(\int_M
dV_{\tilde{g}})^{\frac{n-2}{n}}},$$ where $s_{\tilde{g}}$ and
$dV_{\tilde{g}}$ respectively denote the scalar curvature and the
volume element of $\tilde{g}$. That minimum value is called the
\emph{Yamabe constant} of the conformal class, and denoted as
$Y(M,[g])$. Then the \emph{Yamabe invariant} is defined as the
supremum of the Yamabe constants over the set of all conformal
classes on $M$, and one can hope for a canonical metric as a limit
of such a maximizing sequence.

The Yamabe invariant of a compact orientable surfaces is $4\pi
\chi(M)$ where $\chi(M)$ denotes the Euler characteristic of $M$ by
the Gauss-Bonnet theorem. In general, it is not quite easy to
exactly compute the Yamabe invariant. Recently much progress has
been made in low dimensions. In dimension $3$, the geometrization by
the Ricci flow gave many answers, and in dimension $4$, the
Spin$^{c}$ structure and the Dirac operator have been remarkable
tools for computing the Yamabe invariant. LeBrun \cite{lb1,lb2,lb3}
used the Seberg-Witten theory to show that if $M$ is a compact
K\"ahler surface whose Kodaira dimension $\kappa(M)$ is not equal to
$- \infty$,  then
$$Y(M)=-4\sqrt{2}\pi\sqrt{(2\chi+3\tau)(\tilde{M})},$$ where $\tau$ denotes the
signature and $\tilde{M}$ is the
minimal model of $M$, and for $\Bbb CP^{2}$, $$Y(\Bbb CP^{2})=12\sqrt{2}\pi.$$
In particular, note that if $\kappa(M)=0$ or $1$, $Y(M)=0$.

One notes that the blow-up does not change the Yamabe invariant of
K\"ahler surfaces and may ask:
\begin{question}\label{ques1}
Let $M$ be a smooth closed orientable $4$-manifold with $Y(M)\leq 0$. Is
there an orientation of $M$ such that $Y(M\sharp\ m\ \overline{\Bbb
CP^{2}})=Y(M)$ for any integer $m>0$? What about connected sums or surgeries along circles with $4$-manifolds with negative-definite intersection form and
nonnegative Yamabe invariant?
\end{question}

In this article, we will show :
\begin{thm}\label{th1}
Let $M$ be a closed  K\"ahler surface of $\kappa (M)\geq 0$ (with $b_2^+(M)>1$ if $\kappa (M)>0$), and $N$ be a smooth closed oriented $4$-manifold with  $b_2^+(N)=0$ and $Y(N)\geq 0$.
Then
$$Y(M\# N)=Y(M).$$
\end{thm}
More generally, we prove the case of the surgery along circles.
\begin{defn}
Let $M_1$ and $M_2$ be smooth $n$-manifolds with embedded $k$-spheres $c_1$ and $c_2$ respectively, where the normal bundles are trivial. A {\it surgery} of $M_1$ and $M_2$ along $c_i$'s are defined as the result of deleting tubular neighborhood of each $c_i$ and gluing the remainders by identifying two boundaries $S^k\times S^{n-k-1}$ using a diffeomorphism of $S^k$ and the reflection map of $S^{n-k-1}$. When $M_2$ is not specified, it means a surgery with $S^n$.
\end{defn}

\begin{thm}\label{th3}
Let $M$ be a closed  K\"ahler surface with $\kappa (M)\geq 0$ and $b_2^+(M)>1$, and $N_i$ for $i=1,\cdots,m$ be  smooth closed oriented $4$-manifolds with  $b_2^+(N_i)=0$, $b_1(N_i)\geq 1$, and $Y(N_i)\geq 0$.
Suppose that $c_{i}\subset N_{i}$ is an embedded circle nontrivial in $H_{1}(N_{i},\Bbb R)$ for $i=1,\cdots,m$. If $\tilde{M}$ is a manifold obtained from $M$ by performing a surgery
with $\cup_{i=1}^{m} N_{i}$ along $\cup_{i=1}^{m} c_{i}$, then
$$Y(\tilde{M})=Y(M).$$
\end{thm}
Note that the surgery on $M$ with $(S^1\times S^3)\# N$ along a null-homotopic circle in $M$ and a circle representing $[S^1]\times \{\textrm{pt}\}\in H_1(S^1\times S^3,\Bbb Z)$ is  $M\# N$.

When $b_1(N_i)=0$, theorem \ref{th3} is no longer true in general. For example, take a closed non-spin simply-connected K\"ahler surface $M$ with $\kappa(M)\geq 0$ and $b_2^+(M)>1$. Let $N_i=S^4$, and $c_i$ be any embedded circle in $N_i$ for all $i$. Performing a surgery around null-homotopic circles in $M$ with $\cup_{i=1}^{m} N_{i}$, we get $\tilde{M}$ which is just $M\# (\#_{i=1}^m(S^2\times S^2))$. By applying Wall's results \cite{wall1, wall2}, it is diffeomorphic to $M\# (\#_{i=1}^m(\Bbb CP^2\#\overline{\Bbb CP^2}))$ which again becomes diffeomorphic to $a\ \Bbb CP^{2} \sharp\ a\ \overline{\Bbb CP^{2}}$ where $a=m+\frac{1}{2}(b_{2}(M)+ \tau(M))$, if $m$ is sufficiently large. But $$Y(a\ \Bbb CP^{2} \sharp\ a\ \overline{\Bbb CP^{2}})>0\geq Y(M).$$

We also give a different proof of the following result proved by Gursky and
LeBrun in \cite{Gur}:
\begin{thm}\label{th2}
Let $N$ be a smooth closed oriented $4$-manifold satisfying
$b_2(N)=0$ and $Y(N)\geq
12\sqrt{2}\pi(=Y(\Bbb CP^2))$. Then
$$Y(\Bbb CP^2\# N)=Y(\Bbb CP^2).$$
\end{thm}

For surgeries of codimension less than $3$, in general the Yamabe invariant changes drastically after a surgery. But some surgeries
along $T^{2}$ in $4$-manifolds do preserve the Yamabe invariant.

We introduce some well-known different types of surgeries in $4$-manifolds.
Suppose that a smooth $4$-manifold $M$ contains a homologically essential
tours $T^2$ with self-intersection zero. Deleting a tubular neighborhood
$T^2\times D^2$ of $T^2$ and gluing back using a diffeomorphism $\varphi$ of the boundary $T^3$, we get a new smooth $4$-manifold $M_{\varphi}$ called a generalized logarithmic transform of $M$.

Now suppose that two smooth $4$-manifolds
$M_1$ and $M_2$ each contain an embedded closed surface $F$ with self-intersection zero.
Deleting a tubular neighborhood $F\times D^2$ in each and gluing the remaining parts
along the boundary $F\times S^1$ using a diffeomorphism of $F$ and the complex conjugation map of $S^1$, we get a fiber sum of $M_1$ and $M_2$. When it is performed along two embedded surfaces in $M$, we call it an internal fiber sum of $M$.

Finally a knot surgery manifold for a knot $K\subset S^3$ with the knot exterior $E(K)$ is a smooth $4$-manifold
obtained by gluing $M\setminus (T^2\times D^2)$ and $S^1\times E(K)$ along the boundary $T^3$  in such a way that
the homology class $[\textrm{pt}\times \partial D^2]$ is identified with $[\textrm{pt}\times \lambda]$
where $\lambda$ is a longitude of $K$.
Then a knot surgery of $M$ is the same as the
fiber sum of $M$ with $S^1\times M_K$ along the torus $S^1\times m
\subset S^1\times M_K$, where $m$ is a meridian circle to $K$
and $M_K$ is the $3$-manifold obtained by performing $0$-framed
surgery on $K$.

Now let $M$ be a closed K\"ahler surface of Kodaira dimension equal to
$0$ or $1$ with $b_{2}^{+}(M)>1$. It is known that $M$ admits a
$T$-structure defined by Cheeger and Gromov \cite{CG}. (For an explicit construction, see Paternain and Petean \cite{pp}.) The existence of a $T$-structure implies that the manifold admits a sequence of Riemannian metrics with volume form converging to zero uniformly while the sectional curvatures are bounded below, so that the Yamabe invariant must be nonnegative.(\cite{pp})
Let $\tilde{M}$ be the manifold obtained
from $M$ by a generalized logarithmic transform or an internal fiber sum or a
fiber sum with $S^1\times N$ where $N$ is a closed oriented $3$-manifold with nonzero Seiberg-Witten invariant along an embedded $T^2$ which is a regular orbit of the above $T$-structure. Then $\tilde{M}$ has an obvious induced $T$-structure, and if $\tilde{M}$ also has a nontrivial Seiberg-Witten invariant, we immediately
get $$Y(\tilde{M})=Y(M)=0.$$ It is
interesting to note that these phenomena also appear in some cases
of Kodaira dimension $2$ as follows:
\begin{thm}\label{th4}
Let $M=\Sigma_1\times \Sigma_2$ be a product of two Riemann surfaces
of genus $>1$, and $\alpha_1,\cdots,\alpha_m$ and
$\beta_1,\cdots,\beta_m$ be non-intersecting homologically-essential
circles embedded in $\Sigma_1$ and $\Sigma_2$ respectively.

Suppose that $X_k$ for
$k=1,\cdots, \mu$ where $\mu\leq m$ is a closed oriented $3$-manifold with $b_1(X_k)\geq 1$ and nonzero Seiberg-Witten invariant in a chamber, and $c_k$ for $k=1,\cdots, \mu$ is an embedded circle in $X_k$ representing a non-torsion generator of $H_1(X_k,\Bbb Z)$.

Let
$\tilde{M}$ be a manifold obtained from $M$ by performing  on
$\cup_{i=1}^m \alpha_i \times \beta_i$ an
internal fiber sum or a fiber sum with $\cup_{k=1}^\mu S^1\times X_k$ around $\cup_{k=1}^\mu S^1\times c_k$. Then
$$Y(\tilde{M})=Y(M).$$
\end{thm}

\begin{cor}\label{th5}
Let each $M_i$ for $i=1,\cdots,l$  be a product of two Riemann
surfaces of genus $>1$, and $T_1,\cdots,T_m$ be tori embedded in
$\cup_{i=1}^l M_i$ as above. Let $\tilde{M}$ be a manifold obtained
from $\cup_{i=1}^l M_i$ by performing on
$\cup_{i=1}^m T_i$ an internal fiber sum or
a fiber sum with $S^1\times X_k$'s as above. Then
$$Y(\tilde{M})=-(\sum_{i=1}^l |Y(M_i)|^{2})^{\frac{1}{2}}.$$
\end{cor}
It is left as a further question whether the above theorems still
hold true for any homologically essential tori.

\section{Basic formulae of Yamabe invariant}

When $Y(M)\leq 0$, it can be written as a very nice form:
\begin{eqnarray*}
Y(M)=-\inf_{g}(\int_{M} |s_{g}|^{\frac{n}{2}}
d\mu_{g})^{\frac{2}{n}}=-\inf_{g}(\int_{M} |s_{g}^-|^{\frac{n}{2}}
d\mu_{g})^{\frac{2}{n}},
\end{eqnarray*}
where $s_{g}^-=\min(s_g,0)$. (For a proof, see \cite{lb3,sung3}.)

Another practical formula is the gluing formula of the Yamabe
invariant under the surgery.
\begin{thm}[Kobayashi \cite{koba}, Petean and Yun \cite{PY}]\label{surger}
Let $M_1,M_2$ be smooth closed manifolds of dimension $n\geq 3$.
Suppose that an $(n-q)$-dimensional smooth closed (possibly
disconnected) manifold $W$ embeds into both $M_1$ and $M_2$ with
isomorphic normal bundle. Assume $q \geq 3$. Let $M$ be any
manifold obtained by gluing $M_1$ and $M_2$ along $W$. Then
$$
Y(M) \geq \left\{
\begin{array}{ll}  -( |Y(M_1)|^{n/2}+  |Y(M_2)|^{n/2} )^{2/n}
   &\mbox{if } Y(M_i)\leq 0 \ \forall i\\
  \min(Y(M_1),Y(M_2)) &\mbox{if } Y(M_1) \cdot Y(M_2) \leq 0\\
  \min(Y(M_1),Y(M_2)) &\mbox{if } Y(M_i)\geq 0 \ \forall i\mbox{ and } q=n
\end{array}\right.
$$
\end{thm}

A nontrivial estimation of the Yamabe invariant on  $4$-manifolds
comes from the Seiberg-Witten theory.
\begin{thm}[LeBrun \cite{lb1,lb2}]\label{lebrun}
Let $(M,g)$ be a smooth closed oriented Riemannian $4$-manifold
with $b_2^+(M)\geq 1$. Let $\frak{s}$ be a Spin$^c$ structure on $M$
with first chern class $c_1(\frak{s})$. Suppose that Seiberg-Witten
invariant of $\frak{s}$ is nontrivial in a chamber. Then
\begin{eqnarray*}
Y(M,[g]) \leq  \frac{|4\pi c_1(\frak{s})\cup
[\omega]|}{\sqrt{[\omega]^2/2}}
\end{eqnarray*}
where $\omega$ is nonzero and self-dual harmonic with respect to
$g$. If the Seiberg-Witten invariant of $\frak{s}$ is nontrivial
for any small perturbation, then
\begin{eqnarray*}
Y(M,[g]) \leq -4\sqrt{2}\pi ||c_1^+(\frak{s})||_{L^2}
\end{eqnarray*}
where $c_1^+$ denotes
the self-dual harmonic part of $c_1(\frak{s})$.
\end{thm}

\section{Computation of Seiberg-Witten invariant}
Let $M$ be a smooth closed oriented Riemannian $4$-manifold and $P$ be
its orthonormal frame bundle which is a principal $SO(4)$ bundle.
Consider oriented $\Bbb R^3$-vector bundles $\wedge^2_+$ and
$\wedge^2_-$ consisting of self-dual $2$-forms and
anti-self-dual $2$-forms respectively. Let's let $P_1$ and $P_2$
be associated $SO(3)$ frame bundles. Unless $M$ is spin, it is
impossible to lift these to principal $SU(2)$ bundles. Instead
there always exists the $\Bbb Z_2$-lift, a principal
$U(2)=SU(2)\otimes_{\Bbb Z_2} U(1)$ bundle, of a $SO(3)\oplus
U(1)$ bundle, when the $U(1)$ bundle on the bottom, denoted by $L$,
has first chern class equal to $w_2(TM)$ modulo $2$.  We call this
lifting a Spin$^c$ structure on $M$.

Let $W_+$ and $W_-$ be $\Bbb C^2$-vector bundles associated to
the above-obtained principal $U(2)$ bundles. One can define a
connection $\nabla_A$ on them by lifting the Levi-Civita
connection and a $U(1)$ connection $A$ on $L$. Then the Dirac
operator $D_A : \Gamma(W_+)\rightarrow \Gamma(W_-)$ is defined as
the composition of $\nabla_A:\Gamma(W_+)\rightarrow T^*M\otimes
\Gamma(W_+)$ and the Clifford multiplication.

For a section $\Phi$ of $W_+$,  (perturbed) Seiberg-Witten equations
of $(A,\Phi)$ is given by
$$\left\{
\begin{array}{ll} D_A\Phi=0\\
  F_{A}^+ +\mu= \Phi\otimes\Phi^*-\frac{|\Phi|^2}{2}\textrm{Id},
\end{array}\right.
$$
where $F_A^+$ is the self-dual part of the curvature $dA$ of $A$, and
a purely imaginary self-dual $2$-form $\mu$ is a perturbation term, and finally the identification of both sides in the second equation comes from the Clifford action.

Now we review the Seiberg-Witten invariant as defined by
Ozsv\'ath and Szab\'o \cite{OS}. Suppose $b_2^+(M)>0$, and let $\frak s$
be a Spin$^c$ structure on $M$. The
configuration space $\frak {B}$ of the Seiberg-Witten equations is
given by $$(\mathcal{A}(W_+)\times \Gamma(W_+)) / \textrm{Map}(M,S^1),$$
where $\mathcal{A}(W_+)$ is the space of connections on $L=\det(W_+)$
and is identified with $\Omega^1(M;i\Bbb R)$, and $\textrm{Map}(M,S^1)$ is the group
of gauge transformations. Since $\Gamma(W_+)$ is
contractible, $\frak B$ is homotopy-equivalent to
$T^{b_1(M)}=\frac{H^1(M;\Bbb R)}{H^1(M;\Bbb Z)}$. The irreducible
configuration space $\frak B^*$ is
$$(\mathcal{A}(W_+)\times (\Gamma(W_+)-\{0\}))/ \textrm{Map}(M,S^1),$$
and it is  homotopy-equivalent to $\Bbb CP^\infty \times
\frac{H^1(M;\Bbb R)}{H^1(M;\Bbb Z)}$ so that
$$H^*(\frak B^*;\Bbb Z)\simeq \Bbb Z[U]\otimes \wedge^*H^1(M;\Bbb
Z).$$  Defining the graded algebra $\Bbb A(M)$ over $\Bbb Z$ by
$$\Bbb Z[H_0(M;\Bbb Z)]\otimes \wedge^*H_1(M;\Bbb Z)$$ with
$H_0(M;\Bbb Z)$ grading two and $H_1(M;\Bbb Z)$ grading one, we
have an obvious isomorphism  $$\mu : \Bbb A(M)\tilde{\rightarrow}
H^*(\frak B^*;\Bbb Z)$$ such that $\mu$ maps the positive generator of $H_0(M;\Bbb Z)$
to $U$. Note that the $\mu$ map restricted to a subset $H_1(M;\Bbb
Z)\otimes \Bbb Z$  is given by $Hol_c^*(d\theta)|_{\frak B^*}$ for
$c\in H_1(M;\Bbb Z)$, where $Hol_c: \frak B\rightarrow S^1$ is the
holonomy map around $c$.

Then the Seiberg-Witten invariant
$SW_{M,\frak{s}}$ is a function
$$SW_{M,\frak{s}}:\Bbb A(M)\rightarrow \Bbb Z$$
$$a\mapsto \langle \frak{M}_{M,\frak{s}},\mu(a) \rangle,$$ where
$\frak{M}_{M,\frak{s}}\subset \frak B$ is the moduli space, i.e. the solution space modulo gauge transformations of the Seiberg-Witten equations of $(M,\frak{s})$.
It turns out that $SW_{M,\frak{s}}$ is independent of the Riemannian metric and a generic perturbation, if $b_2^+(M)> 1$. (When $b_2^+(M)=1$, it may depend on the choice of the chamber.) For a noncompact $M$ with cylindrical-end metric, we can do the same job by considering solutions with finite energy. Here, the energy of a solution $(A(t),\Phi(t))$ in temporal gauge on the cylinder $T:=\partial M\times [0,\infty]$ is defined as
$$\| A'(t) \|^2_{L^2(T)}+\| \Phi'(t) \|^2_{L^2(T)},$$
where the temporal gauge means that $A$ has no temporal component $dt$.

We denote $\frak{M}_{M,\frak{s}}^{irr}:=\frak{M}_{M,\frak{s}}\cap \frak B^*$ and $\frak{M}_{M,\frak{s}}^{red}:=\frak{M}_{M,\frak{s}}- \frak{M}_{M,\frak{s}}^{irr}$.
 It is also useful to define the Seiberg-Witten series of $M$ to be the element of the group ring $\Bbb Z[H_2(M,\Bbb Z)]$ given by
$$\overline{SW}_M:=\sum_{\frak s} SW_{M,\frak{s}}(\underbrace{1\otimes \cdots \otimes 1} _{d(\frak s)/2 })\ PD(c_1(\frak s)),$$
where $d(\frak s):=\textrm{dim}_{\Bbb R}\ \frak{M}_{M,\frak{s}}$,  $PD$ denotes  the Poincar\'{e}-dual, and $\frak s$ runs over all isomorphism classes of Spin$^c$ structures on $M$ with even $d(\frak s)$.

For more details about the Seiberg-Witten theory and the gluing of the
moduli spaces, the readers are referred to
\cite{morgan, marcol, mst, nicol1, park, taubes, vid1}.

Before stating the theorem, we note the following lemma.
\begin{lem}
Let $N$ be a smooth closed oriented $4$-manifold with negative intersection form $Q$. Then there exits a Spin$^c$ structure $\frak{s}'$ on $N$ satisfying $c_{1}^{2}(\frak{s}')=-b_2(N)$.
\end{lem}
\begin{proof}
By the Donaldson's theorem, $Q$ is diagonalizable over $\Bbb Z$. (Although the original Donaldson's theorem \cite{donal} is stated for simply-connected ones, a simple Mayer-Vietoris argument can be applied for this generalization.) Let $\{\alpha_1,\cdots,\alpha_{b_2(N)}\}$ be a basis for $H^2(N,\Bbb Z)\otimes \Bbb Q$ diagonalizing $Q$.

We have to show that there exists $x\in H^2(N,\Bbb Z)$ satisfying that
$Q(x,x)=-b_2(N)$, and $x$ is characteristic, i.e. $Q(x,\alpha)\equiv Q(\alpha,\alpha)$ mod 2 for any $\alpha\in H^2(N,\Bbb Z)$. It is easy to check that $x=\sum_{i=1}^{b_2(N)}\pm\alpha_i$ do the job.
\end{proof}

\begin{thm}\label{th3.1}
Let $M$ and $N$ be smooth closed oriented $4$-manifolds such that
$b_2^+(M)>0$, $b_2^+(N)=0$, and $b_1(N)\geq 1$. Let $c\subset N$ be an embedded circle nontrivial in $H_1(N,\Bbb R)$ and $\tilde{M}$ be a manifold obtained by performing a surgery on $M$ with $N$ along $c$.

If $\tilde{\frak{s}}$ is the Spin$^c$
structure on $\tilde{M}$ obtained by gluing a Spin$^c$ structure $\frak{s}$ on $M$ and a Spin$^c$ structure $\frak{s}'$ on $N$ satisfying $c_{1}^{2}(\frak{s}')=-b_2(N)$,
then
$$SW_{\tilde{M},\tilde{\frak{s}}}(a\cdot [d_1]\cdots [d_{b_1(N)-1}])=\pm SW_{M,\frak{s}}(a)$$ for $a\in \Bbb A(M)$, where $[d_1],\cdots, [d_{b_1(N)-1}]$ along with $r[c]$ for some $r\in \Bbb Q$ form a basis for the torsion-free part of  $H_1(N,\Bbb Z)$.
\end{thm}
\begin{proof}
By removing a tubular neighborhood
$S^1\times D^3$ around the circle where the surgery is performed, we construct $\hat{M}$ and $\hat{N}$ with cylindrical end modeled on $S^1\times S^2$ with a standard
metric of positive scalar curvature which we denote by $Y$. For
$Y$ with the trivial Spin$^c$ structure, the moduli space is the set $\chi(Y)$ of flat
connections modulo gauge transformations of the trivial Spin$^c$
structure, which is diffeomorphic to $S^1$.

On $S^1\times D^3$ we put a metric of positive scalar curvature
with the same cylindrical-end, and see that its moduli space with
the trivial Spin$^c$ structure is also the set
$\chi(S^1\times D^3)$ of flat connections modulo gauge
transformations of the trivial Spin$^c$ structure, which is unobstructed.
In an obvious way, $\chi(S^1\times D^3)$ is diffeomorphic to
$\chi(Y)$. From $b_2^+(\hat{M})>0$, $\frak{M}_{\hat{M},\frak{s}}=\frak{M}_{\hat{M},\frak{s}}^{irr}$ and it is unobstructed by using a generic exponentially-decaying perturbation.

Let $\hat{\mathcal{G}}$ be the gauge transformations on $\hat{M}$. (Note that any gauge transformations on $Y$ extend to $S^1\times D^3$ and $\hat{N}$. We will denote such
extensions of $\hat{\mathcal{G}}$ also by $\hat{\mathcal{G}}$ by abuse of notation.) Letting $\hat{\chi}(Y)$ be the set of equivalence classes of flat
connections on $Y$ modulo $\hat{\mathcal{G}}$, $\hat{\chi}(Y)$ is a
covering of $\chi(Y)$ with fiber $H^1(Y,\Bbb Z)/H^1(\hat{M},\Bbb
Z)$. Similarly we define $\hat{\chi}(S^1\times D^3)$ and
$\hat{\frak{M}}_{\hat{N},\frak{s}'}$. (In fact, $\hat{\chi}(S^1\times D^3)=\hat{\chi}(Y).)$  Since the asymptotic map
$$(\partial_\infty,\partial_\infty):
\frak{M}_{\hat{M},\frak{s}}\times \hat{\chi}(S^1\times
D^3)\rightarrow \hat{\chi}(Y)\times \hat{\chi}(Y)$$ is transversal
to the diagonal $\Delta\subset \hat{\chi}(Y)\times \hat{\chi}(Y)$,
$\frak{M}_{M,\frak{s}}$ is diffeomorphic to the fibred product,
i.e
$$\frak{M}_{M,\frak{s}}\simeq(\partial_\infty,\partial_\infty)^{-1}\Delta=
\frak{M}_{\hat{M},\frak{s}}\times_{\hat{\chi}(Y)}
\hat{\chi}(S^1\times D^3)\simeq\frak{M}_{\hat{M},\frak{s}}.$$

For $\hat{N}$ part, first $\hat{\frak{M}}_{\hat{N},\frak{s}'}^{irr}$ is unobstructed by a generic perturbation, but reducible part is nontrivial because $\hat{N}$ does not have a metric of positive scalar curvature in general. Importantly, $\hat{\frak{M}}_{\hat{N},\frak{s}'}^{red}$ is non-empty for any perturbation, because $b_2^+(N)=0$.
\begin{lem}
When $b_1(N)\leq 1$, by a generic exponentially-decaying perturbation,
$\hat{\frak{M}}_{\hat{N},\frak{s}'}^{red}$ is unobstructed for the gluing with $\frak{M}_{\hat{M},\frak{s}}$.
\end{lem}
\begin{proof}
We will follow Vidussi's method \cite{vid2}.
Recall the deformation complex
of appropriate weighted Sobolev spaces :
$$0\rightarrow \Omega_\delta^0(\hat{N},i\Bbb R)\rightarrow
\Omega_\delta^1(\hat{N},i\Bbb R)\times \Gamma_\delta(W_+)\rightarrow
\Omega_{\delta,+}^2(\hat{N},i\Bbb R)\times \Gamma_\delta(W_-)\rightarrow 0$$ and the Kuranishi model near a reducible solution $(A,0)$:
$$H^1(\hat{N},Y;i\Bbb R)\times H^1(Y,i\Bbb R)\times \ker
D_A\rightarrow H^1(Y,i\Bbb R)/H^1(\hat{N},i\Bbb
R)\times\textrm{coker}\ D_A.$$ The virtual dimension of the
moduli space is $$2\ \textrm{ind}_{\Bbb
C}D_A+b_1(\hat{N})=\frac{1}{4}(c_{\hat{N}}-\tau(\hat{N}))-\eta_B(0)+b_1(\hat{N})$$
where $c_{\hat{N}}=-\frac{1}{4\pi^2}\int_{\hat{N}}F_A\wedge F_A$,
$\tau$ is the signature, and $\eta_B(0)$ is the eta invariant of the
Dirac operator associated with the asymptotic limit $B$ of $A$. From
our assumption $c_{\hat{N}}=c_1^2(\frak{s}')=\tau(\hat{N})$, and the
$\eta_B(0)$ vanishes for $Y$ with a standard metric.(see \cite{nicol}.)
Therefore the virtual dimension is $b_1(\hat{N})$. For the surjectivity in the above Kuranishi picture, we only need to show $\textrm{coker}\ D_A=0$ for a generic
exponentially-decaying perturbation. Since the index is zero, it's
equivalent to showing $\ker D_A=0$.

Letting $d^+\nu\in \Omega_{\delta,+}^2(\hat{N},i\Bbb R)$ be a perturbation
term (Recall $b_2^+(\hat{N})=0$.),
$F^+_{A+\nu}=d^+\nu$ and $(A+\nu,0)$ is a reducible solution for the
perturbed Seiberg-Witten equations.  Suppose there exists a nonzero
$\Phi$ satisfying $D_{A+\nu}\Phi=0$. Consider a smooth map
$$F: \frak{M}_{\hat{N},\frak{s}'}^{red}\times (\Gamma_\delta (W_+)-\{0\})\times
\Omega_\delta^1(\hat{N},i\Bbb R)\rightarrow \Gamma_\delta (W_+)$$
$$(A,\Phi,\nu)\mapsto D_{A+\nu}\Phi.$$ Since the differential $DF$
is surjective, $F^{-1}(0)$ is a smooth manifold. Applying the
Sard-Smale theorem to the projection map $\pi_3$ onto the third
factor, for a second category subset of $\nu$, $F^{-1}(0)\cap
\pi_3^{-1}(\nu)$ is a smooth manifold of dimension $b_1(\hat{N})+2\
\textrm{ind}_{\Bbb C} D_{A+\nu}=b_1(\hat{N})\leq 1$. On the other hand, as
$D_{A+\nu}$ is $\Bbb C$-linear, the real dimension of the kernel of
$D_{A+\nu}$ must be greater than or equal to $2$ unless it is empty.
By this contradiction, our claim is proved.
\end{proof}

We first consider the case of $b_1(N)=1$, in which  $\hat{\frak{M}}_{\hat{N},\frak{s}'}^{red}$ is diffeomorphic to
$\hat{\chi}(Y)$, and $\hat{\frak{M}}_{\hat{N},\frak{s}'}^{irr}$ is zero-dimensional by the dimension formula. Chop off $\hat{M}$ and $\hat{N}$ at $Y\times \{t\}$ for $t
\gg 1$ and glue them along the boundary to get $\tilde{M}$. Then
\begin{eqnarray*}
\frak{M}_{\tilde{M},\tilde{\frak{s}}}&\simeq&
(\frak{M}_{\hat{M},\frak{s}}(*)\times_{\hat{\chi}(Y)}(\hat{\frak{M}}_{\hat{N},\frak{s}'}^{irr}
\times S^1))/S^1\cup(\frak{M}_{\hat{M},\frak{s}}\times_{\hat{\chi}(Y)}\hat{\frak{M}}_{\hat{N},\frak{s}'}^{red})\\
&=& (\frak{M}_{\hat{M},\frak{s}}(*)\times_{\hat{\chi}(Y)}\hat{\frak{M}}_{\hat{N},\frak{s}'}^{irr})\cup\frak{M}_{\hat{M},\frak{s}},
\end{eqnarray*}
where $\frak{M}_{\hat{M},\frak{s}}(*)$ is the based moduli space, i.e. the solution space modulo based gauge transformations which are equal to 1 at a fixed point.
As is well-known, the Seiberg-Witten invariant vanishes on $\frak{M}_{\hat{M},\frak{s}}(*)$, because the $\mu$ cocycles are pulled back from $\frak{M}_{\hat{M},\frak{s}}$. Therefore the Seiberg-Witten invariant for $\frak{M}_{\tilde{M},\tilde{\frak{s}}}$ is obtained from the evaluation on $\frak{M}_{\hat{M},\frak{s}}$ to give our desired formula.

Now, we turn to the case when $b_1(N)\geq 2$. Because of the obstruction issue, we first kill $d_i$'s by the surgery, and glue with $M$, and finally  revive the $d_i$'s by the (inverse) surgery. Let $N'$ be the manifold obtained from $N$ by the surgery around $d_1,\cdots,d_{b_1(N)-1}$ with $b_1(N)-1$ copies of $S^4$. Then clearly $b_1(N')=1$ and moreover :
\begin{lem}
Let $U:=N-V$ and $V:=\cup_{i=1}^{b_1(N)-1} S^1\times D^3$ which is a tubular neighborhood of $\cup_{i=1}^{b_1(N)-1} d_i$. Then $H_2(N,\Bbb Z)\simeq H_2(U,\Bbb Z)\simeq H_2(N',\Bbb Z)$ with isomorphic intersection paring, where both isomorphisms are induced by the obvious inclusions.
\end{lem}
\begin{proof}
This can be seen in the Mayer-Vietoris sequence.
First for $(U,V, N=U\cup V)$,
$$H_2(\partial U)\stackrel{i_*}\rightarrow H_2(U)\oplus H_2(V)\stackrel{\varphi}{\rightarrow} H_2(N)\rightarrow
H_1(\partial U)\rightarrow H_1(U)\oplus H_1(V).$$
Since $\varphi$ is surjective, which is  because $H_1(\partial U)$ injects into $H_1(U)$, it is enough to show that $i_*=0$.  Obviously $H_2(V)=0$, and
to prove that $i_*(H_{2}(\partial U))=0\in H_{2}(U)$,
consider the following commutative diagram of exact sequences :
\[
\xymatrix{
H_{3}(U,\partial U)\ar[r]^{\partial_{*}}\ar[d]_{PD} & H_{2}(\partial U)
\ar[r]^{i_{*}}\ar[d]^{PD} & H_{2}(U)\ar[d]^{PD}\\
H^{1}(U)\ar[r]^{i^{*}} & H^{1}(\partial U)\ar[r]^{\partial^{*}} &
H^{2}(U,\partial U).}
\]
Since $c_i$'s are non-torsion in $N$,  $H^{1}(\partial U)=i^*( H^{1}(U))$, and hence it gets mapped to zero by $\partial^{*}$.

The case for $(U,\cup_{i=1}^m D^2\times S^2, N')$ is similar. In this case, $i_*$ maps $H_{2}(\partial U)$ isomorphically onto $H_2(\cup_{i=1}^m D^2\times S^2)$.
\end{proof}

We perform the surgery on $M$ with $N'$ around $c$ to get $M'$. Let $\tilde{\frak{s}}$ be the resulting Spin$^c$ structure on $M'$. (We abused the notation, because it is basically the same as $\tilde{\frak{s}}$ on $\tilde{M}$.) Since $b_2^+(N')=0$, we can apply the previous process to get
$$SW_{M',\tilde{\frak{s}}}(a)= SW_{M,\frak{s}}(a)$$ for $a\in \Bbb A(M)$.

In order to get $\tilde{M}$, we perform an (inverse) surgery on $M'$ around two spheres which are the cores of the added $D^2\times S^2$'s in the surgery around $c_i$'s.
Those two spheres are homologically trivial, and we can apply Ozsv\'ath and Szab\'o's theorem \cite{OS},
$$SW_{\tilde{M},\tilde{\frak{s}}}(a\cdot [d_1]\cdots [d_{b_1(N)-1}])=\pm SW_{M',\tilde{\frak{s}}}(a)$$ for $a\in \Bbb A(M)$. This completes the proof.
\end{proof}

\begin{thm}\label{th3.2}
Let $M$ and $N$ be smooth closed oriented $4$-manifolds such that
$b_2^+(M)>0$ and  $b_2^+(N)=0$. Suppose $\tilde{\frak{s}}$ is the Spin$^c$
structure on $M\# N$ obtained by gluing a Spin$^c$ structure $\frak{s}$ on $M$ and a Spin$^c$ structure $\frak{s}'$ on $N$ satisfying $c_{1}^{2}(\frak{s}')=-b_2(N)$. Then
$$SW_{M\# N,\tilde{\frak{s}}}(a\cdot [d_1]\cdots [d_{b_1(N)}])=\pm SW_{M,\frak{s}}(a)$$ for  $a\in \Bbb A(M)$, where $[d_1],\cdots, [d_{b_1(N)}]$  form a basis for the torsion-free part of  $H_1(N,\Bbb Z)$.
\end{thm}
\begin{proof}
This is an immediate corollary of the previous theorem, because $M\#N$ is the same as the manifold obtained from $M$ by a surgery with $(S^1\times S^3)\# N$ along a circle representing $[S^1]\times \{\textrm{pt}\}\in H_1(S^1\times S^3,\Bbb Z)$.
\end{proof}


\section{Proof of Theorem \ref{th1}}
Just for simplicity, we may assume that $M$ is minimal, because $\overline{\Bbb CP^2}$'s can be absorbed into $N$. By the gluing formula of the
theorem \ref{surger}, $$Y(\tilde{M})\geq Y(M).$$
To obtain the reverse inequality, the computations in the previous
section allows us to apply LeBrun's theorem \ref{lebrun}. Let
$\frak{s}$ be the Spin$^{c}$ structure on $M$ induced by the
canonical line bundle,  which  has nonzero Seiberg-Witten invariant for a chamber.

Let's first consider the case of $\kappa(M)=0$. Recall that $M$ is finitely covered by $T^4$ or K3 surfaces which we denote by $X$.  To the contrary, suppose there exists a metric  of positive scalar curvature on $M\# N$. Then so does $X\# N\# \cdots \# N$ where the number of copies of $N$ is the order of the covering map from $X$ to $M$. Since $b_2^+(X)\geq 2$, the Seiberg-Witten invariant of the obvious Spin$^c$ structure $\tilde{\frak{s}}$ on $X\# N\# \cdots \# N$ is well-defined independently of the chamber and nonzero by theorem \ref{th3.2}. This means that it cannot admit a metric of positive scalar curvature, which is a contradiction.

Now let's consider the case when $\kappa(M)>0$. Let
$c_1(\frak{s})+E$ be the first chern class of $\tilde{\frak{s}}$ on $M\# N$ as in theorem \ref{th3.2},
where $E$  comes from $N$.
For any metric $g$ on $\tilde{M}$
\begin{eqnarray*}
((c_1(\frak{s})\pm E)^{+})^{2}&=&(c_1(\frak{s})^{+}\pm E^{+})^{2} \\
&=&(c_1(\frak{s})^{+})^{2}\pm 2c_1(\frak{s})^{+}\cdot E^{+}+(E^{+})^{2} \\
&\geq& (c_1(\frak{s})^+)^{2}\pm 2c_1(\frak{s})^{+}\cdot E^{+}.
\end{eqnarray*}
Thus at least one of $((c_1(\frak{s})+E)^{+})^{2}$ and
$((c_1(\frak{s})-E)^{+})^{2}$ should be greater than or equal to
$(c_1(\frak{s})^+)^2$.
Say $((c_1(\frak{s})+E)^{+})^{2}\geq  (c_1(\frak{s})^+)^2$.
By applying the second
inequality of the theorem \ref{lebrun}, we get
\begin{eqnarray*}
Y(\tilde{M},[g])&\leq& -4\sqrt{2}\pi ||(c_1(\frak{s})+E)^+||_{L^2}\\ &\leq& -4\sqrt{2}\pi ||c_1^+(\frak{s})||_{L^2}\\ &\leq& -4\sqrt{2}\pi \sqrt{c_1^2(\frak{s})}=Y(M),
\end{eqnarray*}
completing the proof.

\section{Proof of Theorem \ref{th3}}
Again for simplicity's sake, we may assume that $M$ is minimal, because any embedded circle in $X\#\overline{\Bbb CP^2}\# \cdots \# \overline{\Bbb CP^2}$ can be moved to $X-D^4$ by an isotopy, where $D^4$ is the $4$-ball in which the connected sums with $\overline{\Bbb CP^2}$'s are done,  and $\overline{\Bbb CP^2}$'s can be absorbed into $N_1$.
Then the proof is the same as before.

\begin{rmk}
In case that $\kappa(M)=0$ and $b_2^+(M)=1$, if the surgery is done along the circle which is trivially covered by the covering map from $X$ to $M$, then we can lift up the surgery downstairs and use the previous argument in the connected sum case to obtain the same result.
\end{rmk}

\section{Proof of Theorem \ref{th2}}
Again by the gluing formula of the theorem \ref{surger}, it is immediate that
$$Y(\Bbb CP^{2}\# N)\geq Y(\Bbb CP^{2}).$$
For the reverse inequality, let $\frak{s}$ be the Spin$^{c}$
structure on $\Bbb CP^{2}$ induced by the canonical line bundle, and
$[\omega]$ be a nonzero element of $H^{2}(\Bbb CP^{2};\Bbb Z)$. Recall that the
Seiberg-Witten invariant of $(\Bbb CP^{2},\frak{s})$ for a
perturbation $t\omega$ with $|t| \gg 1$ is nonzero for either $t>0$
or $t<0$. By the theorem \ref{th3.2},  so is $(\Bbb CP^{2}\# N,\tilde{\frak{s}})$. Therefore the first inequality of
theorem \ref{lebrun} applies, and the right hand side of the
inequality is
$$\frac{|4\pi c_1\cup
[\omega]|}{\sqrt{[\omega]^2/2}}=\frac{|4\pi (3H\cdot
tH)|}{\sqrt{(tH\cdot tH)/2}}=12 \sqrt{2}\pi,$$ where $H$ denotes the
hyperplane class of $\Bbb CP^{2}$. This completes the proof.

\section{Proof of Theorem \ref{th4} and Corollary \ref{th5}}

Let's first consider the case of the theorem \ref{th4}. Recall that
$M$ admits a K\"ahler-Einstein metric so that
$$Y(M)= -4\sqrt{2}\pi \sqrt{c_1^{2}(\frak{s})},$$ where $\frak{s}$
is the Spin$^c$ structure on $M$ given by the canonical line bundle.
By the adjunction formula, $c_1(\frak{s})$ vanishes on each torus
$T_{j}:=\alpha_j\times \beta_j$.

To apply the product formula of the Seiberg-Witten series, we check
if the so-called "admissibility" condition in \cite{park} is
satisfied. Let's denote $M-(\cup_{j=1}^mT_j\times D^2)$ by $M'$ and
the inclusion map $\partial M' \hookrightarrow M'$ by $i$. Let
$\gamma_j$ be $\{\textrm{pt}\} \times\partial D^2\subset
T_j\times\partial D^2$. There are two non-obvious things to check:
$i_*[\gamma_j]=0 \in H_1(M',\Bbb Z)$  for all $j$, and the
cokernel of $i^*:H^1(M',\Bbb Z)\rightarrow H^1(\partial M',\Bbb Z)$
is freely generated by the Poincar\'e-duals of $[T_j]$'s in $\partial M'$.

For the first one, consider the following commutative diagram of
exact sequences :
\[
\xymatrix{
H_{2}(M',\partial M')\ar[r]^{\partial_{*}}\ar[d]_{PD} & H_{1}(\partial M')
\ar[r]^{i_{*}}\ar[d]^{PD} & H_{1}(M')\ar[d]^{PD}\\
H^{2}(M')\ar[r]^{i^{*}} & H^{2}(\partial M')\ar[r]^{\partial^{*}} &
H^{3}(M',\partial M').}
\]
It's enough to show that
$PD([\gamma_j])$ belongs to the image of $i^*$. This is because
$PD([\gamma_j])\in H^{2}(\partial M')$ which is the dual of
$[T_j]\times \{\textrm{pt}\}\in  H_{2}(\partial M')$ actually comes
from $H^2(M)$ via pull-back.

For the second one, we need the following commutative diagram of
exact sequences :
\[
\xymatrix{
H_{3}(M',\partial M')\ar[r]^{\partial_{*}}\ar[d]_{PD} & H_{2}(\partial M')
\ar[r]^{i_{*}}\ar[d]^{PD} & H_{2}(M')\ar[d]^{PD}\\
H^{1}(M')\ar[r]^{i^{*}} & H^{1}(\partial M')\ar[r]^{\partial^{*}} &
H^{2}(M',\partial M').}
\]
By using the above result $i_*[\gamma_j]=0$, $i_*([\alpha_j]\times
[\gamma_j])$ and $i_*([\beta_j]\times [\gamma_j])$ are all zero in
$H_{2}(M')$. But $i_*([\alpha_j]\times [\beta_j])$ is nonzero
because it is nonzero even in $H_{2}(M)$. Thus the cokernel of $i^*$
is freely generated by $PD([\alpha_j]\times [\beta_j])$'s.

In the same way, these two properties also hold for $X_k':=S^1\times X_k - (S^1\times c_k)$ for all $k$.(Here we need the condition $[c_k]\equiv \pm 1 \in H_1(X_k,\Bbb R)$.)

Note that by using the Mayer-Vietoris argument, it follows from $i_*[\gamma_j]=0$ that $H_2(M')$ is mapped isomorphically into $H_2(M)$ by the inclusion, and likewise for $X_k$'s.

Recall that the Seiberg-Witten series of $M$ is given by
$$\overline{SW}_{M}=[\Sigma_1]^{\chi(\Sigma_2)}[\Sigma_2]^{\chi(\Sigma_1)}+
(-1)^{\frac{\chi(\Sigma_1)\chi(\Sigma_2)}{4}}[\Sigma_1]^{-\chi(\Sigma_2)}[\Sigma_2]^{-\chi(\Sigma_1)}$$
(see \cite{morgan}.),
and the Seiberg-Witten invariant of $S^1\times X_k$ is the same as that of $X_k$ with its basic classes coming from $X_k$ via the pull-back.(see \cite{bal}.)  Now applying the product formula for the Seiberg-Witten series \cite{park},
\begin{eqnarray*}
\overline{SW}_{\tilde{M}}&=&(\overline{SW}_{M'}\prod_{k=1}^\mu
\overline{SW}_{X_k'})|_{\varphi_*}\\
&=&(\overline{SW}_{M}\prod_{j=1}^m
([T_{j}]^{-1}-[T_{j}])\prod_{k=1}^\mu
\overline{SW}_{S^1\times X_k}([S^1\times c_k]^{-1}-[S^1\times c_k]))|_{\varphi_*},
\end{eqnarray*}
where $|_{\varphi_*}$ denotes the identification in the homology induced by the gluing map of the fiber sum construction, and  if $b_1(X_k)=1$, we mean the Seiberg-Witten Series for a chamber.


Now taking $\frak{s}'$ with nonzero Seiberg-Witten invariant from $S^1\times X_k$ parts  and gluing with $\frak s$, we obtain a Spin$^c$ structure $\tilde{\frak{s}}$ with nonzero Seiberg-Witten invariant on $\tilde{M}$ such that  $c_{1}^{2}(\tilde{ \frak{s}})=c_{1}^{2}(\frak{s})$, because $$c_{1}^{2}(\frak{s'})=\langle c_1(\frak s),c_1(\frak{s}')\rangle = \langle c_1(\frak s),T_j\rangle= \langle c_1(\frak s'),T_j\rangle=[T_j]\cdot [T_k]=0\ \ \forall j,k.$$ This enables us to apply the second inequality of the theorem \ref{lebrun} and to get
$$Y(\tilde{M})\leq  -4\sqrt{2}\pi \sqrt{c_1^2(\frak{s})}=Y(M).$$


To show the reverse inequality, we need to construct a Riemannian
metric on $\tilde{M}$ whose Yamabe constant is arbitrarily close to
$Y(M)$. Let's take a maximal subset of
$\{\alpha_1,\cdots,\alpha_{m}\}$, any two elements of which are mutually non-isotopic,
and may assume that it is $\{\alpha_1,\cdots,\alpha_{m'}\}$ for
$m'\leq m$ by renaming. In the same way, we define
$\{\beta_1,\cdots,\beta_{m''}\}$. Let $g_{1}$ be a complete metric
of constant curvature $-1$ on $\hat{\Sigma}_{1}:=\Sigma_{1}-\cup_{j=1}^{m'}\alpha_j$.
It is well-known that the metric near the infinity is the cusp
metric, i.e. $dt^2+e^{-2t}g_{S^1}, t\in [a,\infty)$, where $g_{S^1}$ is the metric on
the circle of radius $1$. At each cusp,
we cut it at $t=b$ for $b\gg 1$ and  glue a cylinder with a metric
$dt^2+ e^{-2b}g_{S^1}, t\in [b,b+1]$ along $\{b\}\times S^1$. Then
the resulting metric is only $C^0$, so to obtain a nearby smooth
metric, take a smooth decreasing convex function
 $\rho:[b-1,b]\rightarrow [0,1]$ such that $\rho\equiv e^{-t}$ near
$b-1$, and $\rho\equiv e^{-b}$ near $b$. Then $dt^2+\rho^2 g_{S^1}$
is a smooth metric with curvature ranging from $-1$ to $0$, and we glue the
corresponding cylindrical ends along the boundary to get back
$\Sigma_1$ with a metric $\tilde{g_1}$ parameterized by $b\gg 1$. In
the same fashion, we construct $\tilde{g_2}$ on $\Sigma_2$
parameterized by $b\gg 1$, using a complete metric $g_2$
of constant curvature $-1$ on $\hat{\Sigma}_{2}:=\Sigma_{2}-\cup_{j=1}^{m''}\beta_j$.

In $(M,\tilde{g_1}+\tilde{g_2})$, we can find a
$\delta$-neighborhood $N_j=\{x\in M | \textrm{dist}(x,T_j)\leq
\delta\}$ for all $j=1,\cdots,m$ such that they are mutually
disjoint for some $\delta>0$ when $b$ and $c$ are sufficiently
large. Note that $N_j$ are all isometric to the product
$e^{-2b}g_{S^1}+e^{-2b}g_{S^1}+ g_{D^2(\delta)}$ where
$g_{D^2(\delta)}$ is the flat metric on the disk of radius $\delta$,
and $\delta$ can remain constant if we take $b$ further
larger. For the fiber sum, we bend $g_{D^2(\delta)}$ on $D^2(\frac{\delta}{2})-D^2(\frac{\delta}{3})$ to a metric like a horn with a
cylindrical end $dt^2+g_{S^1(\frac{\delta}{3})}, t\in [0,\frac{\delta}{4}]$, where $g_{S^1(\frac{\delta}{3})}$ is the metric on the
circle of radius $\frac{\delta}{3}$.

For the $S^1\times X_k-(S^1\times c_k)$ part, we first consider a finite-volume complete metric $h$ on
$X_k- c_k$ such that $h$ near the infinity is of the form $e^{-2t}g_{S^1}+dt^2+g_{S^1(\frac{\delta}{3})}$, where $e^{-2t}g_{S^1}$ is the metric in the $c_k$ direction.  Then applying the cutoff procedure as before, we change $h$ into a metric $\tilde{h}$ parameterized by $b$ with a cylindrical end $e^{-2b}g_{S^1}+dt^2+g_{S^1(\frac{\delta}{3})},  t\in [0,\frac{\delta}{4}]$. Surely the volume and curvature of $\tilde{h}$ is bounded independently of $b>0$. We finally take the metric $e^{-2b}g_{S^1}+\tilde{h}$ on $S^1\times (X_k-(S^1\times c_k))$.

We now perform the fiber sum to get a metric $\tilde{g}$ on
$\tilde{M}$. The important thing is that if we take $b$ sufficiently large, the volume
of the gluing region and the parts from $S^1\times X_k$'s is made arbitrarily small with its curvature bounded.
Thus applying the Gauss-Bonnet theorem for complete finite-volume hyperbolic
surfaces, we have that for any $\epsilon>0$, there exists $\tilde{g}$ such that
\begin{eqnarray*}
-(\int_{\tilde{M}} s_{\tilde{g}}^2
d\mu_{\tilde{g}})^{\frac{1}{2}}&\geq&-(\int_{\hat{\Sigma}_{2}\times \hat{\Sigma}_{2}} s_{g_1+g_2}^{2}
d\mu_{g_1+g_2})^{\frac{1}{2}}-\epsilon \\ &=&
-2(4\pi\chi(\hat{\Sigma}_{1})
4\pi\chi(\hat{\Sigma}_{2}))^{\frac{1}{2}}-\epsilon
\\ &=&
-2(4\pi\chi(\Sigma_{1}) 4\pi\chi(\Sigma_{1}))^{\frac{1}{2}}-\epsilon
\\ &=& Y(M)-\epsilon,
\end{eqnarray*}
which is our desired inequality.


The case of Corollary \ref{th5} goes exactly the same. What we need is
Kobayashi's formula \cite{koba} on the Yamabe invariant of the
disjoint union by which
$$Y(M_1\cup \cdots \cup M_l)=-(\sum_{i=1}^l
|Y(M_i)|^{2})^{\frac{1}{2}}$$ for $Y(M_i)\leq 0$ $\ \forall i$.

\begin{rmk}
As mentioned in the introduction, a knot surgery is a special case of the above
construction. $M_K$ has the same homology as $S^1\times S^2$ with $[m]$ generating $H_1(M_K,\Bbb Z)$, and $$\overline{SW}_{S^1\times M_K}=\frac{\Delta_{K}([T]^{2})}{([T]^{-1}-[T])^2},$$ where $T$ denotes $S^1\times m$, and  $\Delta_{K}$ is the symmetrized Alexander polynomial of a knot $K$. (For a proof, see \cite{FS} and \cite{park}.)

It is an interesting question whether a knot surgery on $M=\Sigma_1\times \Sigma_2$ as above does not change the homeomorphism class of $M$.
\end{rmk}

\section{Examples}
Let $M$ be a K\"ahler surface of nonnegative Kodaira dimension, and
$N_{i}$  be an $S^{1}$ bundle over a rational homology $3$-sphere
for $i=1,\cdots,m$. Then $$Y(M\# N_1 \# \cdots \# N_m)=Y(M).$$ Also
we can perform  surgeries with a product of $S^1$ with a rational
homology $3$-sphere along $S^1\times \{\textrm{pt}\}$ to get the
same result.

For $\Bbb CP^{2}$ case, presently we don't have many examples but
$$Y(\Bbb CP^{2}\#_{i=1}^{m} (S^{1}\times S^{3}))=12 \sqrt{2}\pi.$$

\bigskip

\noindent{\bf Acknowledgement.} The author warmly thanks Dr.
Young-Heon Kim for helpful discussions on the Teichm\"uller space, and all NIMS people.


\begin{thebibliography}{99}



\bibitem{bal} S. Baldridge, {\em Seiberg-Witten invariants, orbifolds, and circle actions}, Trans. of AMS, {\bf 355} (2002), 1669--1697.

\bibitem{besse} A. Besse, Einstein Manifolds, Springer-Verlag (1987).

\bibitem{CG} J. Cheeger and M. Gromov, {\em Collapsing Riemannian manifolds
while keeping their curvature bounded I},
J. Diff. Geom. {\bf 23} (1986), 309--346.

\bibitem{donal} S. Donaldson, {\em An application of gauge theory to four dimensional topology}, J. Diff. Geom. {\bf 18} (1983), 279--315.

\bibitem{FS} R. Fintushel and R. Stern,
{\em  Knots, links, and $4$-manifolds}, Invent. Math. {\bf 134} (1998), 363--400.


\bibitem{Gur} M.J. Gursky and C. LeBrun,
{\em Yamabe constants and spin$^c$ structures}, GAFA {\bf 8} (1998), 965--977.

\bibitem{koba} O. Kobayashi, {\em Scalar curvature of a metric with unit
volume}, Math. Ann. {\bf 279} (1987), 253--265.

\bibitem{lb1} C. LeBrun,
{\em Four manifolds without Einstein metrics}, Math. Res. Lett.
{\bf 3} (1996), 133--147.

\bibitem{lb2} C LeBrun, {\em Yamabe constants and the perturbed Seiberg-Witten
equations}, Comm. Anal. Geom.  {\bf 5} (1997), 535--553.

\bibitem{lb3} C LeBrun, {\em Kodaira dimension and the Yamabe problem},
Comm. Anal. Geom.  {\bf 7} (1999), 133--156.


\bibitem{marcol} M. Marcolli, {\em Seiberg-Witten-Floer homology and Heegaard splittings}, Internat. J. Math. {\bf 7} (1996), 671--696.

\bibitem{morgan} J. Morgan, {\em The Seiberg-Witten Equations and
applications to the topology of smooth four-manifolds}, Princeton
University Press, 1996.

\bibitem{mst} J. W. Morgan, Z. Szabo, and C. H. Taubes
{\em A product formula for the Seiberg--Witten invariants and the generalized Thom conjecture},  J. Diff. Geom. {\bf 44} (1996), 706--788.


\bibitem{nicol} L.I. Nicolaescu, {\em Eta invariants of Dirac
operators on circle bundles over Rieman surfaces and virtual
dimensions of finite energy Seiberg-Witten moduli spaces}, Israel
J. Math. {\bf 114} (1999), 61--123.

\bibitem{nicol1} L.I. Nicolaescu, Notes on Seiberg-Witten theory, AMS (2000).

\bibitem{park} B. Doug Park, {\em A Gluing formula for the
Seiberg-Witten invariant along $T^3$},  Michigan Math. J. {\bf 50}
(2002), 593--612.

\bibitem{pp} G. P. Paternain and J. Petean, {\em Minimal entropy and
collapsing with curvature bounded from below}, Invent. Math. {\bf
151} (2003), 415--450.

\bibitem{PY} J. Petean and G. Yun, {\em Surgery and the Yamabe invariant},
GAFA {\bf 9} (1999), 1189--1199.

\bibitem{OS} P. Ozsv\'ath and Z. Szab\'o, {\em Higher type adjunction
inequalities in Seiberg-Witten theory},  J. Diff. Geom.  {\bf 55}
(2000), 385--440.

\bibitem{sung1} C. Sung, {\em Surgery, curvature, and minimal volume},
Ann. Global Anal. Geom. {\bf 26} (2004), 209--229.

\bibitem{sung2} C. Sung, {\em Surgery and equivariant Yamabe invariant},
Diff. Geom. Appl. {\bf 24} (2006), 271--287.

\bibitem{sung3} C. Sung, {\em Collapsing and monopole classes of $3$-manifolds},
J. Geom. Phys. {\bf 57} (2007), 549--559.

\bibitem{taubes} C. H. Taubes, {\em  The Seiberg--Witten invariants and $4$-manifolds
with essential tori}, Geom. Topol. {\bf 5} (2001), 441--519.

\bibitem{vid1} S. Vidussi, {\em Seiberg-Witten theory for $4$-manifolds decomposed
along  $3$-manifolds of positive scalar curvature},
Pr\'epublication \'Ecole Polytechnique {\bf 99-5} (1999).

\bibitem{vid2} S. Vidussi, {\em Seiberg-Witten invariants for manifolds diffeomorphic
outside a circle}, Proc. A.M.S. {\bf 129} (2001), 2489--2496.

\bibitem{wall1} C. T. C. Wall, {\em Diffeomorphisms of $4$-manifolds},
J. London Math. {\bf 39} (1964), 131--140.

\bibitem{wall2} C. T. C. Wall, {\em On simply-connected $4$-manifolds},
J. London Math. {\bf 39} (1964), 141--149.
\end{thebibliography}
\end{document}